\documentclass[12pt]{article}

\def\Endproof{{\ \vbox{\hrule\hbox{%
   \vrule height1.3ex\hskip0.8ex\vrule}\hrule
}}}
\def\0{\leqno}
\author{{\small by}\\Gabriel Mitric }
\date{}
\title{POISSON STRUCTURES ON COTANGENT BUNDLES
}
\begin{document}
\maketitle
{\def\thefootnote{*}\footnotetext[1]%
{{\it 2000 Mathematics Subject Classification}: 53D17.
Cotangent bundle. Horizontal lift. }}
\begin{center} \begin{minipage}{12cm} A{\footnotesize BSTRACT.
 We make a study of Poisson structures of
$T^{*}M$ which are graded structures when  restricted   to the
fiberwise polynomial algebra, and give examples.
 A class of more general graded bivector fields
which induce a given Poisson structure $w$  on the base
 manifold $M$ is constructed.
 In particular, the {\it horizontal lifting} of a Poisson structure
  from $M$ to $T^{*}M$ via  connections gives  such bivector
  fields and we discuss the conditions for these lifts to be Poisson
  bivector fields and their compatibility with the canonical
  Poisson structure on $T^{*}M.$
  Finally, for a $2-$form $\omega $ on a Riemannian manifold,
 we study the conditions for some associated $2-$forms
 of $\omega $ on  $T^{*}M$
 to define Poisson structures on cotangent bundles.
  }
\end{minipage}
\end{center} \noindent

\section*
{\begin{center} {\bf 1. Graded Poisson structures on cotangent
bundles}
\end{center}}
Let $M$ be an $n-$dimensional differentiable manifold and
 $\pi :T^{*}M\longrightarrow M$ its cotangent bundle.
If  $(x^i),$  $(i=1,...,n),$ are local coordinates on $M,$ we
denote by $(p_i)$ the covector coordinates with respect to the
cobasis $(dx^i).$ (We assume that everything is $C^{\infty}$ in
this paper).

In this section we discuss  {\it graded} Poisson structures $W$ on
the cotangent bundle $T^{*}M$
 obtained as {\it lifts} of Poisson structures $w$ on the base manifold $M$,
 in the sense that  the canonical projection $\pi $
 is a Poisson mapping (see \cite{GMIV}).

Denote by $S_k(TM)$ the space of $k-$contravariant symmetric
tensor fields on $M$ and by $\odot$ the symmetric tensor product
on the algebra  $S(TM)=\bigoplus \limits _{k\geq 0}S_k(TM).$ The
spaces of fiberwise homogeneous $k-$polynomials
$$
{\mathcal{H}\mathcal{P}}_k(T^{*}M):=\{\tilde{Q}=Q^{i_1...i_k}p_{i_1}...
p_{i_k} /
  \0(1.1)
$$
$$
  Q=Q^{i_1...i_k}\frac{\partial}{\partial x^{i_1}}\odot ... \odot
\frac{\partial}{\partial x^{i_k}} \in
 S_k(TM)\} \
$$
are interesting subspaces of the function space
$C^{\infty}(T^{*}M),$ and will play an important role in this
paper. ($:=$ denotes a definition).

The map
$$
\sim \ : (S(TM),\odot )\longrightarrow ( {\mathcal{P}}(T^{*}M), \
\cdot \ ) \ , \ \ \sim Q:=\tilde{Q} \ ,
  \0(1.2)
 $$
where $ {\mathcal{P}}(T^{*}M):={\oplus
}_k{\mathcal{H}}{\mathcal{P}}_k(T^{*}M)$ is the {\it polynomial
algebra} and the dot denotes the usual multiplication, is an
isomorphism of algebras.

 On $T^{*}M$ we also have the spaces of (fiberwise) non homogeneous
 polynomials of degree $\leq k$
 $$
{\mathcal{P}}_{k}(T^{*}M):= \bigoplus \limits
_{h=0}^k{\mathcal{H}\mathcal{P}}_{h} \ .
 $$

For $k=1,$ ${\mathcal{A}}(T^{*}M):={\mathcal{P}}_{1}(T^{*}M)$ is
the space  of {\it affine functions}, having the elements of the
form
$$
a(x,p)=f(x)+m(X) \ ,
$$
where $f\in C^{\infty}(M),$ $X\in \chi (M)$ (the space of vector
fields on $M$) and $m(X):=\sim X$ is the {\it momentum} of $X.$
($m(X)$ is $X$ regarded as a function on $T^{*}M$).

The elements of the space ${\mathcal{P}}_2(T^{*}M)$
 of non-homogeneous quadratic polynomials are
$$
t(x,p)=f(x)+m(X )+s(Q) \ ,
$$
where $Q=Q^{ij}({\partial}/{\partial x^i}) \odot
 ({\partial}/{\partial x^j})$ is a symmetric contravariant tensor
field on $M$ and $s(Q):=\sim Q.$

Hereafter, by a polynomial on $T^{*}M$ we always mean a fiberwise
polynomial. Also, we write $f$ for both $f$ on $M,$ and $f\circ
\pi $ on $T^{*}M.$

 DEFINITION 1.1. A Poisson structure $W$ on $T^{*}M$ is called {\it
polynomially graded} if $\forall  Q,R\in {\mathcal{P}}(T^{*}M)$
$$
Q\in {\mathcal{P}}_{h}, \ R\in {\mathcal{P}}_k \Longrightarrow
\{Q,R\}_W\in {\mathcal{P}}_{h+k} \ .
  \0(1.3)
$$

PROPOSITION 1.2. {\it A polynomially graded Poisson structure $W$
on $T^{*}M$ induces a Poisson structure $w$
 on the base manifold $M,$ such that the projection
 $\pi :(T^{*}M,W)\longrightarrow (M,w)$ is a Poisson mapping.   }

{\it Proof.} Any function   $f$ on $M$ is a polynomial $(f\circ
\pi ) \in {\mathcal{P}}_{0}(T^{*}M).$   By $(1.3),$ $\forall f,g
\in C^{\infty}(M),$ $\{f\circ \pi, g\circ \pi \}_W\in
C^{\infty}(M)$ and
$$
\{f,g\}_w:=\{f\circ \pi, g\circ \pi \}_W,
     \0(1.4)
$$
defines a Poisson structure $w$ on $M.$ \Endproof

Hereafter, the bracket $\{ \ , \ \}_W$ will be denoted  simply by
$\{ \ , \ \}.$

 If the local coordinate expression of the Poisson structure $w$
 introduced by Proposition 1.2   is
 $$
w=\frac {1}{2}w^{ij}(x)\frac {\partial}{\partial x^i}\wedge \frac
{\partial}{\partial x^j} \ ,
   \0(1.5)
 $$
 Definition 1.1 tells us that $W$ must have the local coordinate
expression
$$
W=\frac {1}{2}w^{ij}(x)\frac {\partial}{\partial x^i}\wedge \frac
{\partial}{\partial x^j}+ (\varphi ^{i}_j(x)+p_aA^{ia}_j(x))\frac
{\partial}{\partial x^i}\wedge \frac {\partial}{\partial p_j}+
    \0(1.6)
$$
$$
+\frac {1}{2}(\eta _{ij}(x)+p_aB
^{a}_{ij}(x)+p_ap_bC^{ab}_{ij}(x))\frac {\partial}{\partial
p_i}\wedge \frac {\partial}{\partial p_j} \ ,
$$
where $w,$ $\varphi ,$  $\eta ,$ $A ,$ $B,$ $C$  are local
functions on $M.$

The Poisson structure $W$ is completely determined by the brackets
$\{ f,g \},$ $\{ m(X),f \}$ and $\{m(X),m(Y)  \},$ where $f,g\in
C^{\infty}(M)$ and $X,Y \in \chi (M),$ since the local coordinates
$x^i$ and $p_i$ are functions of this type $(p_i =m(
{\partial}/{\partial x^i})).$

By $(1.3),$ the bracket $\{m(X),f\}$ is in
${\mathcal{P}}_1(T^{*}M),$ i.e.,
$$
\{ m(X),f \}=Z_Xf+m({\gamma}_Xf) \ ,
   \0(1.7)
$$
where $Z_Xf\in C^{\infty}(M)$ and ${\gamma}_Xf \in \chi (M).$

 $\{m(X),\ . \ \}$ is a derivation of
$C^{\infty}(M). $ Hence,
 $ Z_X$
is  a vector field on $M,$ and the mapping $
  \gamma _{X}:C^{\infty}(M)\longrightarrow \chi (M)$
  also  is a derivation. Therefore, $\gamma _{X}f $ depends only on $df.$

From the Leibniz rule we get that $Z_{hX}=hZ_X$ $(h\in
C^{\infty}(M))$ and $\gamma $ must satisfy
$$
\gamma _{hX}f=h\gamma _{X}f+ (X_h^{w}f)X \ .
   \0(1.8)
$$

The bracket of two affine functions has an expression of the form
$$
\{m(X),m(Y) \}=\beta (X,Y)+ m(V(X,Y))+s(\Psi (X,Y)) \ ,
   \0(1.9)
$$
where $\beta (X,Y) \in C^{\infty}(M),$ $V(X,Y)\in \chi (M)$ and
$\Psi (X,Y) \in S_2(TM)$ are skew-symmetric operators. If we
replace $Y$ by $fY$ in $(1.9),$ the Leibniz rule gives that $\beta
$ is a $2-$form on $M$ and
$$\begin{array}{l}
V(X,fY)=fV(X,Y)+(Z_Xf)Y, \vspace{2mm}\\
   \Psi (X,fY)=f\Psi
(X,Y)+({\gamma}_Xf)\odot Y. \end{array}
  \0(1.10)
$$

DEFINITION 1.3. A  polynomially graded Poisson structure $W$ on
$T^{*}M$ is said to be a {\it graded structure} if $\forall
 Q\in {\mathcal{H}\mathcal{P}}_{h},$ $\forall
 R\in {\mathcal{H}\mathcal{P}}_{k},$
 $\{Q,R\}_W\in {\mathcal{H}\mathcal{P}}_{h+k} \ .$

Remark that a polynomially graded structure on $T^{*}M$ is graded
iff $Z_X=0,$ $\beta =0$ and $V=0.$ In this case $(1.6)$ reduces to
$$
W=\frac {1}{2}w^{ij}(x)\frac {\partial}{\partial x^i}\wedge \frac
{\partial}{\partial x^j}+ p_aA^{ia}_j(x)\frac {\partial}{\partial
x^i}\wedge \frac {\partial}{\partial p_j}
    \0(1.11)
$$
$$
+\frac {1}{2}p_ap_bC^{ab}_{ij}(x)\frac {\partial}{\partial
p_i}\wedge \frac {\partial}{\partial p_j} \ .
$$
As in \cite{GMIV},  a bivector field $W$ on $T^{*}M$ which is
locally of the form $(1.6)$ (respectively $(1.11)$) is called a
{\it polynomially graded} (respectively {\it graded) bivector
field.}

PROPOSITION 1.4.  {\it If $W$ is a graded bivector field on
$T^{*}M$ which is $\pi -$related with a Poisson structure $w$ on
$M,$ there exists a contravariant connection $D$ on the Poisson
manifold $(M,w)$ such that
$$
\{m(X),f  \}=-m(D_{df}X) \ , \ \ \ X\in \chi (M), \ f\in
C^{\infty}(M) \ .
    \0(1.12)
$$
Moreover, if $W$ is a graded Poisson structure on $T^{*}M$ then
the connection $D$ is flat. }

{\it Proof.} A contravariant connection on $(M,w)$ is a
contravariant derivative on  $TM$ with respect to the Poisson
structure \cite{Va1}.

The required connection is defined by
$$
D_{df}X:=-{\gamma}_Xf \ .
   \0(1.13)
$$

That we really get a connection, which is flat in the Poisson
case, follows in exactly the same way as in \cite{GMIV}.
${\Endproof }$

The relation $(1.12)$ extends to

 PROPOSITION 1.5. {\it If $Q$ is a symmetric contravariant tensor
field on $M$ and $\tilde{Q}$ is its corresponding polynomial then,
for any graded Poisson bivector field $W$ on $T^{*}M$, one has  }
$$
\{ \tilde{Q},f\}_W=-\widetilde{D_{df}Q} \ .
  \0(1.14)
$$

{\it Proof.}  $D_{df}$ of $(1.14)$  is  extended  to $S(TM)$ by
$$
(D_{df}Q)(\alpha _1,...,\alpha _k)=X_f^{w}(Q(\alpha _1,...,\alpha
_k))-\sum \limits _{i=1}^kQ(\alpha _1,...,D_{df}\alpha
_i,...,\alpha _k) \ ,
$$
where $\alpha _1,...,\alpha _k\in {\Omega }^1 (M),$ and
$D_{df}\alpha $ is defined by
$$
<D_{df}\alpha ,X >=X_f^{w}<\alpha ,X>-<\alpha ,D_{df}X>, \ \ \
X\in \chi (M) \ .
$$
We put
$$
D_{dx^i}\frac{\partial}{\partial x^j}=-\Gamma
^{ik}_j\frac{\partial}{\partial x^k} \ ,
   \0(1.15)
$$
and by a straightforward computation we get for $\{\tilde{Q},f\}$
 and $-\widetilde{(D_{df}Q)}$ the same local coordinate
 expression. (See \cite{GMIV} for the complete proof in the case
 of a symmetric covariant tensor field on $M.$)
 \Endproof

In order to discuss the next two  Jacobi identities,  let us make
some remarks concerning the operator $\Psi $
 of $(1.9),$ which is given in the case of a
  graded Poisson structure on $T^{*}M$ by
$$
\{m(X),m(Y )\}= s(\Psi (X ,Y ))  \ , \ \ \ X,Y \in \chi (M) \ .
  \0(1.16)
$$
With $(1.13),$ the second relation   $(1.10)$ becomes
$$
\Psi (X , fY )=f\Psi (X , Y ) -\frac {1}{2} (D_{df}X \otimes Y + Y
\otimes D_{df}X  ) \
 \0(1.17)
$$
and this allows us to derive the local coordinate expression of
$\Psi .$ If $X=X^i({\partial}/{\partial x^i })$ and
 $Y=Y^j({\partial}/{\partial x^j }),$ we obtain
$$
\Psi (X,Y)=X^iY^j\Psi \left( \frac{\partial}{\partial x^i },
\frac{\partial}{\partial x^j } \right)
   \0(1.18)
$$
$$
+\left( X^h\frac{\partial Y^j}{\partial x^k }\Gamma ^{ki}_h-
 Y^h\frac{\partial X^i}{\partial x^k }\Gamma ^{kj}_h \right)
\frac{\partial }{\partial x^i}\odot \frac{\partial }{\partial x^j
} +w^{kh}\frac{\partial X^i}{\partial x^k } \frac{\partial
Y^j}{\partial x^h }\frac{\partial }{\partial x^i}\odot
\frac{\partial }{\partial x^j} \ .
$$
Remark that   $\Psi :TM\times TM\longrightarrow {\odot}^2TM $ is a
bidifferential operator of the first order.

 PROPOSITION 1.6. {\it If we define an operator $D_{df}$
which acts on   $\Psi $ by
$$
(D_{df}\Psi )(X ,Y ):=D_{df}(\Psi (X ,Y))- \Psi (D_{df}X , Y)-\Psi
(X ,D_{df} Y) \ ,
   \0(1.19)
$$
the  Jacobi identity
$$
\{\{m(X),m(Y)\},f\}+\{\{m(Y),f\},m(X)\} + \{\{ f,m(X)\},m(Y)\}=0
  \0(1.20)
$$
has the  equivalent form }
$$
(D_{df}\Psi )(X ,Y) =0 \ , \ \ \forall X ,Y \in \chi (M) \ .
  \0(1.21)
$$

{\it Proof. } Using $(1.12),$ $(1.14)$ and $(1.16)$ for $Q=\Psi
(X,Y), $ $(1.20)$ becomes $(1.21).$
 \Endproof

We also find
$$(D_{df}\Psi ) (X ,h Y)=h(D_{df}\Psi )(X , Y)-
[C_D(df,dh)X ]\odot Y \
  \0(1.22)
$$
and hence, we see that  $(1.21)$  is invariant by $X \mapsto fX, \
$ $Y \mapsto gY $ ($f,g\in C^{\infty}(M)$) iff the curvature
$C_D=0.$

Concerning the Jacobi identity
$$
\sum \limits _{(X ,Y ,Z  )} \{ \{m (X),m(Y) \}, m(Z) \} =0 \ ,
 \0(1.23)
$$
(putting indices between parentheses denotes that summation is on
cyclic permutations of these indices)  remark that one must have
an operator $ \Theta $ such that
$$
\{ s(G), m (X) \}=\widetilde{\Theta (G, X)} \ , \ \ X\in \chi (M),
\ G\in S_2(M) \ ,
  \0(1.24)
$$
and
 $\Theta (G, X)$ is a symmetric
$3-$contravariant tensor field  on $M.$

We get the formula
$$
\Theta (fG,hX  )=fh\Theta (G,X)-f(D_{dh}G  )\odot X +hG\odot
D_{df}X+ \{ f,h \}_{w}G \odot X  ,
   \0(1.25)
$$
and then, the local coordinate expression
$$
\Theta (G,X)=G^{ij}X^k\Theta \left(\frac{\partial}{\partial
x^i}\odot \frac{\partial}{\partial x^j},
 \frac{\partial}{\partial x^k}\right)+\frac {1}{3}\sum \limits _{(i,j,k)}
(G^{hj}\frac{\partial X^k}{\partial x^a}\Gamma ^{ai}_h
    \0(1.26)
$$
$$
+G^{ih}\frac{\partial X^k}{\partial x^a}\Gamma ^{aj}_h -
\frac{\partial G^{ij}}{\partial x^a}X^h\Gamma ^{ak}_h+ w^{ab}
\frac{\partial G^{ij}}{\partial x^a} \frac{\partial X^k}{\partial
x^b})\frac{\partial}{\partial x^i}\odot \frac{\partial}{\partial
x^j} \odot \frac{\partial}{\partial x^k} \ .
$$

 Using the operator $\Theta ,$ the Jacobi identity $(1.23)$ becomes
$$
\sum \limits _{(X ,Y , Z  )}\Theta (\Psi (X ,Y), Z )=0 \ ,
 \0(1.27)
$$
and we may summarize our analysis  concerning the graded Poisson
structures on $T^{*}M$ in

 PROPOSITION 1.7. {\it A graded Poisson structure $W$ on $T^{*}M$
 with the bracket $\{ \ , \  \}$ is defined by

 a) a Poisson structure $w$ on the base manifold $M$ such that
 $$
\{ f,g \}_{W}=\{ f,g \}_w \ , \ \ \ f,g \in C^{\infty}(M) \ .
 $$

b) a flat contravariant connection $D$ on $(M,w)$ such that
$$
\{  m(X),f \}=-m( D_{df}X) \ , \ \ X\in C^{\infty}(M) \ .
$$

c) an operator $\Psi :TM\times TM\longrightarrow {\odot}^2TM $
such that
$$
\{m(X),m(Y )\}= s(\Psi (X ,Y ))  \ , \ \ \ X,Y \in \chi (M) \
$$
and formula $(1.21)$ holds.

d) an operator $\Theta $ defined by $(1.24),$ satisfying
$(1.27).$}

\vspace{3mm} To give examples, we consider the following
situation, similar
 to \cite{GMIV}.

Let $(M,w)$ be an $n-$dimensional  Poisson manifold and  suppose
that its symplectic foliation $S$ is contained in a regular
foliation ${\mathcal{F}}$ on $M$, such that $T{\mathcal{F}}$ is a
{\it foliated bundle}
 i.e.,  there are local bases
$\{Y_u\}$ $(u=1,...,p ,$ $p= rank \ \mathcal{F})$
 of $T\mathcal{F}$ with  transition  functions constant along
 the leaves of $\mathcal{F}.$
 Consider a decomposition
$$
TM=T{\mathcal{F}}\oplus \nu  \mathcal{F} \ ,
 \0(1.28)
 $$
where $\nu  \mathcal{F}$ is a complementary subbundle of $T
\mathcal{F},$ and  $\mathcal{F}-$adapted local coordinates
$(x^a,y^u) \ \ (a=1,...,n-p) $
 on $M$ \cite{Va2}.

The Poisson bivector $w$  has the form
 $$
w=\frac {1}{2}w^{uv}(x,y)\frac {\partial}{\partial y^u}\wedge
\frac {\partial}{\partial y^v} \ \ \ \ \ \ (w^{vu}=-w^{uv}) \ ,
  \0(1.29)
 $$
since $S\subseteq {\mathcal{F}}.$

If $\{ \beta ^u \},$ $\{ {\tilde{\beta}}^v\}$ $(u,v=1,...,p)$
 are the dual cobases of $\{Y_u \},$ $\{
{\tilde{Y}}_v \}$ $(\beta ^u (Y_v)=\delta ^u_v)$ then  their
transition functions are constant along the leaves of
${\mathcal{F}}.$

Now, $\forall \alpha \in T^{*}M,$ $\alpha =\zeta
_adx^a+\varepsilon _u\beta ^u$ and we may consider $(x^a,y^u,\zeta
_a, \varepsilon _u)$ as {\it distinguished local coordinates } on
$T^{*}M.$ The transition function are
$$
{\tilde{x}}^a ={\tilde{x}}^a(x), \
{\tilde{y}}^u={\tilde{y}}^u(x,y), \ {\tilde{\zeta}}_u
=\frac{\partial x^a}{\partial {\tilde{x}}^u}\zeta _a \ , \
 \ {\tilde{\varepsilon}}_u =a^v_u(x)\varepsilon _v \ .
  \0(1.30)
$$

PROPOSITION 1.8. {\it Under  the previous hypotheses, $W$ given
with respect to the distinguished local coordinates by
$$
W=\frac {1}{2}w^{uv}(x,y)\frac {\partial}{\partial y^u}\wedge
\frac {\partial}{\partial y^v} \
     \0(1.31)
$$
defines  a  graded Poisson bivector on $T^{*}M.$}

 {\it Proof.} From $(1.30)$ it follows   that
$W$ of $(1.31)$ is a global tensor field on $T^{*}M.$ The
Schouten-Nijenhuis bracket $[W,W]$ has the same expression as
$[w,w]$ on $M$ and thus, the Poisson condition $[W,W]=0$ holds.

To prove that $W$ is graded, we also consider natural coordinates
and show that the
expression of $W$ with respect to these coordinates becomes of the
form $(1.11)$ (see \cite{GMIV}). \Endproof

 There are some interesting particular cases of
Proposition 1.8:

a) $w$ is a regular Poisson structure, and the bundle $TS$ is a
foliated bundle; in this case  we may  take $ {\mathcal{F}}=S.$

b) $S$ is contained in a regular foliation ${\mathcal{F}}$ which
admits adapted local coordinates $(x^a,y^u)$ with local transition
functions
$$
{\tilde{y}}^v=p^v_u(x)y^u+q^v(x) \ .
$$
(${\mathcal{F}}$ is a
 leaf-wise, locally affine, regular foliation.)
In this case $({\partial}/{\partial y^u})=\sum \limits _v a^v_u(x)
({\partial}/{\partial {\tilde{y}}^v})$ and we may use the local
vector fields $Y_u= {\partial}/{\partial y^u} \ .$

c) There exists a flat linear connection $\nabla $ (possibly with
torsion) on the Poisson manifold $(M,w).$ In this case we may
consider as leaves of $ {\mathcal{F}}$ the connected components of
$M,$ and the local $\nabla -$parallel vector fields have constant
transition functions along these leaves. Therefore, we may take
them as $Y_i$ $(i=1,...,n).$

In particular, we have  the result of c)  for a locally affine
manifold $M$ (where $\nabla $ has no torsion), using as $Y_i$
local $\nabla -$parallel vector fields,  and also  for a
parallelizable manifold $M$ (where we have global vector fields
$Y_i$).

As a consequence,  Proposition 2.8 holds for the Lie-Poisson
structure \cite{Va1} of any  dual  ${\mathcal{G}}^{*}$ of a Lie
algebra ${\mathcal{G}},$ the graded Poisson structure
 being defined on  $T^{*}{\mathcal{G}}^{*}={\mathcal{G}}^{*} \times
{\mathcal{G}}.$

\section*
{\begin{center}{\bf 2. Graded bivector fields on cotangent
bundles}\end{center}} In this section we will discuss
 graded  bivector fields on a cotangent bundle $T^{*}M,$
 which may be seen as
lifts of a given Poisson structure $w$ on $M,$ that satisfy less
restrictive existence conditions than in the case of graded
Poisson structures.

Recall  the following definition from \cite{GMIV}. Let
${\mathcal{F}}$ be an arbitrary regular foliation, with
$p-$dimensional leaves, on an $n-$dimensional manifold $N.$ We
denote by $C^{\infty}_{fol}(N)$ the  space of {\it foliated
functions} (the functions on $N$ which are constant along the
leaves
 of $\mathcal{F}$).
 A {\it transversal Poisson  structure} of
$(N,\mathcal{F})$ is  a bivector field  $w$ on $N$  such that
$$
\{f,g\}:=w(df,dg)  \ , \ f,g\in C^{\infty}_{fol}(N)
   \0(2.1)
$$
is a Lie algebra bracket on  $C^{\infty}_{fol}(N).$  A bivector
field  $w$ on $N$   defines a
   transversal Poisson  structure  of $(N,\mathcal{F})$
   iff \cite{GMIV}
 $$
  ({\mathcal{L}}_Yw)|_{ Ann
\ T{\mathcal{F}}}=0 \ ,  \ \ \
 [w,w]|_{ Ann \ T{\mathcal{F}}}=0 \ ,
   \0(2.2)
$$
for all $Y\in \Gamma (T{\mathcal{F}})$ (the space of global cross
sections of $T{\mathcal{F}}),$ where
 $Ann \ T{\mathcal{F}}\subseteq \Omega ^1(N)$ is the
annihilator space  of $T\mathcal{F}.$ ($\Omega ^1(N)$ denotes the
space of Pfaff forms on $N.$)

The cotangent bundle $T^{*}M$ of any manifold $M$  has the
vertical foliation $\mathcal{F}$ by fibers with the tangent
distribution $V:=T\mathcal{F}.$

Obviously, the set of  foliated functions on $T^{*}M$ may be
identified with $C^{\infty}(M).$

PROPOSITION 2.1. {\it  Any polynomially  graded bivector field
 $W$ on $T^{*}M,$ which is $\pi $ related with a Poisson structure of $M$ is
 a transversal Poisson structure of $(T^{*}M,V).$    }

{\it Proof.} The local coordinate expression of $W$ is of the form
$(1.6),$ and $W$ is  $\pi -$related with the bivector field  $w$
defined on $M$ by the first term of $(1.6).$ Then, $(2.2)$ holds,
because $w$ is a Poisson bivector on $M.$ \Endproof

DEFINITION 2.2. A transversal Poisson  structure of
 the vertical foliation  of $T^{*}M$ will be called a
{\it semi-Poisson structure} on $T^{*}M.$

REMARK 2.3.  The structures $W$  of Proposition 2.1 are
polynomially graded semi-Poisson structures on $T^{*}M.$

In what follows, we will discuss some interesting classes of
graded semi-Poisson  structures of $T^{*}M.$ Then, we give a
method to  construct all the  graded semi-Poisson bivector fields
 on $T^{*}M$ which  induce the same  Poisson
 structure $w$ on the base manifold $M.$

Let $D$ be a contravariant derivative on a Poisson manifold
$(M,w).$ First, $\forall Q\in S_k(TM),$ define $^s {D } Q\in
S_{k+1}(TM)$ by
$$
(^s D Q) (\alpha _1,...,\alpha _{k+1})=\frac {1}{k+1}\sum \limits
_{i=1}^{k+1} (D _{\alpha _i}Q)(\alpha _1,...,\hat{\alpha
_i}...\alpha _{k+1}) \ ,
     \0(2.3)
$$
where $\alpha _1,...,\alpha _{k+1} \in \Omega ^1(M),$ and the hat
denotes the absence of the corresponding factor.

If $X=X^i({\partial}/{\partial x^i})\in \chi (M)$  then $DX,$
defined by $(DX)({\alpha }_1,{\alpha }_2)=(D_{\alpha _1}X)\alpha
_2,$  is a $2-$contravariant tensor field on $M,$ and
$$
DX=D^iX^j \frac{\partial}{\partial x^i}\otimes
    \frac{\partial}{\partial x^j} \ ,
      \0(2.4)
$$
where $D^iX^j=(D_{dx^i}X)dx^j=D_{dx^i}X^j -X(D_{dx^i}dx^j).$
According to $(1.15)$ we must have
$$
D_{dx^i}dx^j=\Gamma ^{ij}_k dx^k \ ,
    \0(2.5)
$$
and obtain
$$
D^iX^j={(dx^i)}^{\sharp}X^j-\Gamma ^{ij}_kdx^k=\{
x^i,X^j\}_{w}-\Gamma ^{ij}_kX^k \ .
    \0(2.6)
$$
 Then
$$
^s D X= \frac {1}{2}(D^iX^j+D^jX^i)\frac{\partial}{\partial
x^i}\odot
    \frac{\partial}{\partial x^j}
$$
 and we get
$$
^s D X= \frac {1}{2}[\{ x^i,X^j\}_{w}+\{x^j, X^i \}_{w}-\Gamma
^{ij}_kX^k-\Gamma ^{ji}_kX^k] \frac{\partial}{\partial x^i}\odot
    \frac{\partial}{\partial x^j} \ .
    \0(2.7)
$$

PROPOSITION 2.4. {\it Let $(M,w)$ be a Poisson manifold and $D$ a
contravariant derivative of $(M,w).$ The bivector field   $W_1$ on
$T^{*}M,$
 of bracket $ \{ \ , \    \}_{W_1}$  defined by the conditions
$$
\{f,g\}_{W_1}:=\{f,g\}_{w} \ ,
 \0(2.8)
$$
$$
\{m(X),f\}_{W_1}:=-m(D_{df}X ) \ ,
 \0(2.9)
$$
$$
\{m(X),m(Y)\}_{W_1}=\frac 12 s[^{s} {D } < X,Y>
     \0(2.10)
$$
$$-<^{s} {D}
X,Y>-<X ,^{s} {D} Y >] \ ,
$$
where $f,g\in C^{\infty}(M), \ $ $X,Y \in \chi (M)$ and $< \ , \
>$ is the Schouten-Nijenhuis bracket of symmetric tensor fields
(defined by the natural Lie algebroid of $M$) \cite{{Ba},{GMIV}}
defines a graded semi-Poisson structure on $T^{*}M$ which is $\pi
-$related with $w.$ }

{\it Proof.} If the local coordinate expression of $w$ is $(1.5),$
using $(2.7)$ and the properties of $< , >$ \cite{{Ba},{GMIV}} we
get
$$
W_1=\frac {1}{2}w^{ij}\frac {\partial}{\partial x^i}\wedge \frac
{\partial}{\partial x^j}-p_a \Gamma ^{ia}_{j}\frac
{\partial}{\partial x^i}\wedge \frac {\partial}{\partial p_j}
    \0(2.11)
$$
$$
-\frac 14 p_ap_b\left[\frac {\partial}{\partial x^j}(\Gamma
^{ab}_i+\Gamma ^{ba}_i)-\frac {\partial}{\partial x^i}(\Gamma
^{ab}_j+\Gamma ^{ba}_j)\right] \frac {\partial}{\partial
p_i}\wedge \frac {\partial}{\partial p_j} \ . \  { \Endproof}
$$

REMARK 2.5. The relation $(2.10)$ provides us the expression of
the operator $\Psi _{W_1}$ associated to $W_1$ (see $(1.16)$):
$$
{\Psi}_{W_1}(X,Y )= \frac 12(  ^{s} {D } < X,Y>-<^{s} {D} X,Y>-<X
,^{s} {D} Y >   ) \ .
 \0(2.12)
$$

\vspace{2mm} Now, instead of $D$ we consider a linear connection
$\nabla $ on a Poisson manifold $(M,w)$ and define the vector
field $K$ on $T^{*}M$ by
$$
K (\alpha )={({\sharp }_w \alpha )}^H_{\alpha } \ , \ \ \ \alpha
\in T^{*}M \ ,
     \0(2.13)
$$
where  $\sharp _{w}:T^{*}M\longrightarrow TM$ is defined by
 $ \beta({\alpha}^{\sharp } )=w(\alpha , \beta ), \ $ $\forall \beta
\in {\Omega}^1(M), $
and
the upper index $H$ denotes the horizontal lift with respect
to $\nabla $ (\cite{{Leon},{KY}}). In local coordinates we get
$$
K=p_aw^{ai}\frac {\partial}{\partial x^i}
+\frac{1}{2}p_ap_b(w^{ak}\Gamma ^b_{ki}+w^{bk}\Gamma ^a_{ki})
\frac {\partial}{\partial p_i} \ .
     \0(2.14)
$$

On $T^{*}M$ we have  the canonical symplectic form $ \omega
=d\lambda =dp_i\wedge dx^i \ $ where $\lambda =p_idx^i$ is the
Liouville form, and  the vector bundle isomorphism
$$
{\sharp}_{\omega}:T^{*}M\longrightarrow TM \ , \ \ \ i_X\omega \in
T^{*}M \longmapsto X\in TM
$$
leads to the canonical  Poisson bivector $ W_0:={\sharp}_{\omega}
\omega $ on $T^{*}M.$ It follows
$$
W_0(dF,dG)=\omega (\sharp (dF),\sharp (dG))) \ , \ \ F,G\in
C^{\infty}(T^{*}M) \ ,
    \0(2.15)
$$
and locally one has
$$
W_0=\frac {\partial}{\partial p_i}\wedge \frac {\partial}{\partial
x^i} \ .
  \0(2.16)
$$

PROPOSITION 2.6. {\it If $(M,w)$ is a Poisson manifold then the
bivector field
$$
W_2=\frac 12{\mathcal{L}}_KW_0 \ ,
   \0(2.17)
$$
defines a graded semi-Poisson  structure on $T^{*}M$ which is $\pi
-$related with $w.$ }

{\it Proof.}  We get
$$
W_2=  \frac {1}{2}w^{ij}\frac {\partial}{\partial x^i}\wedge \frac
{\partial}{\partial x^j} +\frac 12p_a({\nabla
}_jw^{ai}+2w^{ik}\Gamma ^a_{kj})
\frac {\partial}{\partial x^i}\wedge \frac {\partial}{\partial
p_j}
    \0(2.18)
$$
$$
+\frac 14 p_ap_b\left[\frac {\partial}{\partial x^j}(w^{ak}\Gamma
^{b}_{ki}+w^{bk}\Gamma ^{a}_{ki})-\frac {\partial}{\partial
x^i}(w^{ak}\Gamma ^{b}_{kj}+w^{bk}\Gamma ^{a}_{kj})\right] \frac
{\partial}{\partial p_i}\wedge \frac {\partial}{\partial p_j} \ ,
$$
where ${\nabla }_jw^{ai}$ are the components of the $(2,1)-$tensor
field on $M$ defined by $X\mapsto {\nabla }_Xw, \ $ $X\in \chi
(M).$ \Endproof

We will say that $W_2$ of $(2.17)$ is the {\it graded $\nabla
-$lift} of the Poisson structure $w$ of $M.$

Using local coordinates and the notation of $(1.2)$ we get
$$
{\mathcal{L}}_K\tilde{Q}= \widetilde{^s D Q  } \ ,
   \0(2.19)
$$
where $D$ is the contravariant derivative induced by the linear
connection $\nabla ,$  defined by $D_{df}={\nabla
}_{{(df)}^{\sharp}}$ (see \cite{Va1}).

From $(2.17)$ we have
$$
\{F_1,F_2\}_{W_2}:=W_2(dF_1,dF_2)
 \0(2.20)
$$
$$
=\frac {1}{2}({\mathcal{L}}_K(\{F_1,F_2\}_{W_0})
-\{{\mathcal{L}}_KF_1,F_2\}_{W_0}-\{F_1,{\mathcal{L}}_KF_2\}_{W_0})
\ , \
$$
where $F_1,F_2 \in C^{\infty}(T^{*}M).$

If $Q_1,Q_2 \in S(TM),$ using $(2.19)$ and the relation
$$
\{ \tilde{Q},\tilde{H}\}_{W_0}:=\widetilde{<Q,H>} \ , \ \ \ Q,H
\in S(TM)
     \0(2.21)
$$
(see  \cite{{Ba},{GMIV}}) we get the explicit formula
$$
\{{\tilde{Q}}_1, {\tilde{Q}}_2   \}_{W_2}=\frac 12\sim[^{s} {D
}<Q_1,Q_2>
    \0(2.22)
$$
$$
-<^{s} {D } Q_1,Q_2>-<Q_1,^{s} {D} Q_2>] \  .
$$

PROPOSITION 2.7. {\it The graded $\nabla -$lift  $W_2$ of $w$
 is characterized by:

$i)$ the Poisson structure induced on $M$ by $W_2$  is  $w,$
 i.e.
$$
\{f,g\}_{W_2}=\{f,g\}_{w} \ , \  \ \forall f,g \in C^{\infty}(M);
 \0(2.23)
$$

$ii)$ for every $f \in C^{\infty}(M)$ and $X \in \chi (M)$
$$
\{m(X),f\}_{W_2}=-m({\bar{D}}_{df}X ) \ ,
 \0(2.24)
$$
where $\bar{D}$ is the contravariant derivative of $(M,w)$ defined
by
$$
{\bar{D}}_{\alpha}\beta=D_{\alpha }\beta+\frac
12({\nabla}_{\cdot}w)(\alpha ,\beta ) \ , \ \ \alpha, \beta \in
\Omega ^1(M) \ ,
   \0(2.25)
$$
 the contravariant derivative $D$ is induced by $\nabla $
 and $({\nabla}_{\cdot}w)(\alpha ,\beta ) $
is the $1-$form $X\mapsto ({\nabla }_Xw)(\alpha ,\beta).$

$iii)$ for any vector fields $X$ and $Y$  of $M$ we have
$$
\{m(X),m(Y)\}_{W_2}=\frac 12 s(^{s} {D} < X,Y>
 \0(2.26)
$$
$$
-<^{s} {D } X ,Y >-<X ,^{s} {D } Y >) \ .
$$
}

{\it Proof.} $i)$ If $f\in C^{\infty}(M)$ then $Df=-X_f^w$ and
from $(2.20),$ $(2.21)$ and the formula
$$
<Q,f>=i(df)Q \ , \ \ f\in C^{\infty}(M), \ Q\in S_p(TM) \ ,
$$
 we get
$$
\{ f,g\}_{W_2}=-\frac 12(<Df,g>+<f,Dg>)=\frac
12(X_f^wg-X_g^wf)=\{f,g \}_w \ .
$$

$ii)$ As $W_2$ is graded, the bracket $\{ m(X),f \}_{W_2}$ must be
of the form $(2.24).$ Denoting
$$
{\bar{D}}_{dx^i}dx^j={\bar{\Gamma}}^{ij}_kdx^k \
$$
$(2.18)$ give us
$$
{\bar{\Gamma}}^{ij}_k= {\Gamma} ^{ij}_k +\frac 12{\nabla}_kw^{ij}
\ ,
     \0(2.27)
$$
where
$$
{\Gamma }^{ij}_k=-w^{ih}{\Gamma}^j_{hk} \ ,
     \0(2.28)
$$
(${\Gamma}^i_{jk}$ are the coefficients of the linear connection
$\nabla $) and hence $(2.25).$

$iii)$ $(2.26)$ is a direct consequence of $(2.22).$ \Endproof

Notice from $(2.26)$ that  the operator ${\Psi}_{W_2}$ associated
to $W_2$ has the same expression as ${\Psi}_{W_2}$ of $(2.12),$
but in the case of $W_1$ the contravariant derivative $D$ is
induced by a linear connection $\nabla $ on $M.$

PROPOSITION 2.8. {\it If the graded semi-Poisson structure $W_1$
is defined by a linear connection on $(M,w)$ then it coincides
with $W_2$ iff $w$ is $\nabla -$parallel. }

{\it Proof.} Compare the characteristic conditions of Propositions
2.4 and 2.7 (or  the coefficients of $({\partial}/{\partial
x^i})\wedge ( {\partial}/{\partial p_j})$
 of $(2.11)$ and of $(2.18),$ using $(2.28)$).
\Endproof

We will prove now

PROPOSITION 2.9. {\it Let $(M,w)$ be a Poisson manifold and
 $\pi :T^{*}M \longrightarrow M$ its cotangent bundle.
 The graded semi-Poisson structures $W$ on $T^{*}M$  which
are $\pi -$related with $w$  are defined by the relations
$$
\{f,g\}_W=\{f,g\}_{w} \ , \ \{m(X),f\}_W=-m(D_{df}X ) \ \mbox{and}
 \
$$
$$
\{m(X),m(Y)\}_W =s(\Psi (X ,Y )) \ , \ \ \ f,g \in C^{\infty}(M),
\  X,Y \in \chi (M) \
$$
 where $D$ is an
arbitrary contravariant connection of $(M,w)$ and the operator
$\Psi $ is given by
$$
\Psi ={\Psi}_0 +A+T \ ,
 \0(2.29)
$$
where ${\Psi}_0$ is the
  operator $\Psi $ of a fixed  graded semi-Poisson structure,
 $A:TM\times TM\longrightarrow {\odot}^2TM$ is a
skew-symmetric, first order, bidifferential operator such that
$$
A(X ,fY )=f A(X , Y)-{\tau}(df ,X )\odot Y \ ,
 \0(2.30)
$$
where $\tau $ is a  $(2,1)-$tensor field on $M, $ and $T$ is a
$(2,2)-$tensor field on $M$ with the properties $T(Y,X)=-T(X,Y)$
and $T(X,Y)\in S_2(TM), \ $ $\forall X,Y \in \chi (M).$ }

 {\it Proof.} If two graded semi-Poisson bivector fields,
$\pi -$related with $w,$ have associated the same contravariant
connection $D,$ it follows from $(1.17)$ that the difference
${\Psi}^{'}-\Psi $ is a tensor field $T,$ as in Proposition. To
change $D$ means to pass to a  contravariant connection
$D^{'}=D+\tau ,$ where $\tau $ is a $(2,1)-$tensor field on $M$
and  from $(1.17)$ again, it follows that $A={\Psi}^{'}-\Psi $
becomes a bidifferential operator with the property $(2.29).$
\Endproof

\section*
{\begin{center}{\bf 3. Horizontal lifts of  Poisson
structures}\end{center}} In this section we define and study an
interesting class of semi-Poisson structures on $T^{*}M$ which are
produced by a process of {\it horizontal lifting } of Poisson
structures from $M$ to $T^{*}M$ via  connections.

 On $T^{*}M$ we  distinguish the vertical distribution $\mathrm{V},$
 tangent to the fibers of the projection $\pi $ and, by complementing
 $\mathrm{V}$ by a  distribution $\mathrm{H},$ called
 {\it horizontal,}
we define a {\it nonlinear connection} on $T^{*}M$ \cite{{Vo1},
{Vo2}}.

We have ({\it adapted}) bases of the form
$$
\mathrm{V}=span \left\{  \frac {\partial}{\partial p_i}  \right\}
\ , \ \mathrm{H}=span \left\{  \frac {\delta}{\delta x^i}=\frac
{\partial}{\partial x^i}-N_{ij}\frac {\partial}{\partial p_j}
\right\} \ ,
   \0(3.1)
$$
and $N_{ij}$ are the  {\it coefficients of the connection} defined
by $\mathrm{H}.$

Equivalently, a nonlinear connection may be seen as an almost
product structure $\Gamma $ on $T^{*}M$ such that
 the eigendistribution   corresponding to the eigenvalue $-1$ is the vertical
 distribution ${\mathrm{V}}$ \cite{Vo1}.

We assume that  the nonlinear connection above is symmetric, i.e.,
$N_{ji}=N_{ij}.$ This condition is independent \cite{Vo1} on the
local coordinates.

The complete integrability of  $\mathrm{H},$ in the sense of the
Frobenius theorem, is equivalent to the vanishing of the curvature
tensor field
 $$
R=R_{kij}dx^i\wedge dx^j \otimes \frac {\partial}{\partial p_k} \
, \ \ R_{kij}=\frac {\delta N_{kj}}{\delta x^i}- \frac {\delta
N_{ki}}{\delta x^j} \ .
    \0(3.2)
 $$
 For a later utilization, we also notice  the formulas
 \cite{{Vo1},{Vo2}}
$$
\left[\frac {\delta}{\delta x^i},\frac {\delta}{\delta
x^j}\right]=-R_{kij} \frac {\partial}{\partial p_k} \ , \ \
\left[\frac {\delta}{\delta x^i},\frac {\partial}{\partial
p_j}\right]=-\Phi ^j_{ik} \frac {\partial}{\partial p_k} \ ,
   \ \Phi ^j_{ik}=-\frac {\partial N_{ik}}{\partial
p_j} \ .
 \0(3.3)
$$

Let $w$ be a bivector on  $M,$ with the local coordinate
expression $(1.5).$

 DEFINITION 3.1.
The {\it horizontal lift} of $w$  to the cotangent bundle $T^{*}M$
is the (global) bivector field $w^H$ defined by
$$
w^H=\frac 12 w^{ij}(x)\frac {\delta}{\delta x^i}\wedge \frac
{\delta}{\delta x^j} \ .
 \0(3.4)
$$

PROPOSITION 3.2. {\it Let $(M,w)$ be a Poisson manifold. If the
connection  $\Gamma$ on $T^{*}M$ is defined by  a linear
connection $\nabla $ on $M,$  the bivector $w^H$  defines a
 graded semi-Poisson structure on $T^{*}M.$}

 {\it Proof.}
  In this case the coefficients of $\Gamma$
 are
 $$
 N_{ij}=-p_k{\Gamma}^k_{ij},
     \0(3.5)
 $$
  where ${\Gamma }^k_{ij}$ are
 the coefficients of $\nabla $ and, with respect to the bases
 $\{ {\partial}/{\partial x^i},
  {\partial}/{\partial p_j}\},$ the local expression
  of $w^H$ becomes
$$
W=\frac {1}{2}w^{ij}\frac {\partial}{\partial x^i}\wedge \frac
{\partial}{\partial x^j}+ w^{ik}\Gamma ^a_{kj}p_a\frac
{\partial}{\partial x^i}\wedge \frac {\partial}{\partial p_j}
 \0(3.6)
$$
$$
+\frac 12 w^{kh}\Gamma ^a_{ki}\Gamma ^b_{hj}p_ap_b \frac
{\partial}{\partial p_i}\wedge \frac {\partial}{\partial p_j} \ .
\     {\Endproof}
$$

 PROPOSITION 3.3. {\it The horizontal lift
 $w^H$ is a Poisson bivector on the
cotangent bundle $T^{*}M$ iff $w$ is a Poisson bivector on the
base manifold $M$ and
$$
R(X_f^H,X_g^H)=0, \ \  \ \forall f,g \in C^{\infty }(M) \ ,
   \0(3.7)
$$
where  $X_f^H$ denotes the usual horizontal lift
\cite{{Leon},{KY}}, from $M$ to $T^{*}M,$ of the $w-$Hamil\-tonian
vector field $X_f$ on $M.$

In this case,  the projection $\pi :(T^{*}M,w^H) \longrightarrow
(M,w)$ is a Poisson mapping. }

 {\it Proof.} We compute the bracket $[w^H,w^H]$ with respect to
 the bases $(3.1)$ and get that the Poisson condition
 $[w^H,w^H]=0$ is equivalent with the pair of conditions
 $$
\sum \limits _{(i,j,k)}w^{hk} \frac {\partial w^{ij}}{\partial
x^h}=0 \ , \ \  w^{il}w^{jh}R_{klh}=0 \ .
     \0(3.8)
 $$
(Putting indices between parentheses denotes that summation is on
cyclic permutations of these indices.)
\newline
The first condition  $(3.8)$ is equivalent to $[w,w]=0$ and the
second is the local coordinate expression of  $(3.7).$ \Endproof

Notice that the condition $(3.7)$ has the equivalent form
$$
R({(\sharp \alpha)}^H,{(\sharp \beta)}^H)=0, \ \ \forall \alpha,
\beta \in {\Omega}^1(M) \ .
 \0(3.9)
$$

REMARK 3.4.  If $w$ is defined by a symplectic form on $M,$
condition $(3.8)$ becomes $R=0.$

 COROLLARY 3.5. {\it If $(M,w)$ is a Poisson manifold and
the connection $\Gamma $ on $T^{*}M$ is defined by a linear
connection $\nabla $ on $M,$ the bivector $w^H$  defines a
  Poisson structure on $T^{*}M$ iff the curvature
  $C_D$ of the contravariant
  connection induced by $\nabla$ on $TM$ vanishes.}

{\it Proof.} If $R^h_{kij}$ are the components of the curvature $R_{\nabla }$
then
$$
R_{kij}=-p_hR^h_{kij} \
    \0(3.10)
$$
and $(3.9)$  becomes
$$
R_{\nabla}(\sharp \alpha ,\sharp \beta)Z=0 \ , \ \ \forall \alpha,
\beta \in {\Omega}^1(M) \ , \ \forall Z \in \chi (M) \ ,
 \0(3.11)
$$
(or, equivalently
$$
R_{\nabla}(X_f,X_g)Z=0 \ , \ \ \forall f,g\in C^{\infty}(M), \
\forall Z \in \chi (M) \ .)
 \0(3.11')
$$
This is equivalent to $C_D=0.$ \Endproof

In the case where $w^H$ is a Poisson bivector, it is interesting
to study its compatibility with the canonical Poisson structure
$W_0$ of $(2.15).$

  PROPOSITION 3.6. {\it If $w^H$  is a Poisson
 bivector, then it is compatible with $W_0$ iff }
 $$
\frac {\partial w^{ij}}{\partial x^k}+w^{ih}\Phi ^j_{hk}-
w^{jh}\Phi ^i_{hk} =0 \ , \ \ w^{ih}R_{hjk}=0 \ .
 \0(3.12)
 $$

{\it Proof.} By a straightforward  computation we get that
 the compatibility condition $[w^H,W]=0$ is equivalent to $(3.12).$
 \Endproof

The Bianchi identity \cite{Vo1}
$$
R_{kij}+R_{ijk}+R_{jki}=0 \ ,
    \0(3.13)
$$
shows that the second relation  $(3.12)$ implies $(3.7).$ Then

  COROLLARY 3.7. {\it If $(M,w)$ is a Poisson manifold and the
cotangent bundle $T^{*}M$ is endowed with a symmetric nonlinear
connection, then $w^H$ is a Poisson bivector on $T^{*}M$
compatible with $W_0$ iff  conditions $(3.12)$  hold.}

 REMARK 3.8. Considering the isomorphism
$$
\Psi :{\mathrm{V}}_u\longrightarrow {\mathrm{H}}_u^{*} \ , \ \
\Psi (X_k {\partial}/{\partial p_k})=X_kdq^k \ ,
$$
where $u\in T^{*}M$ and ${\mathrm{H}}_u^{*}$ is the dual space of
${\mathrm{H}}_u,$ the second condition $(3.12)$ may be written in
the equivalent form
$$
[\Psi(R(X,Y))]{({\sharp}_w\alpha)}^H=0 \ , \ \ \forall X,Y \in
\chi (T^{*}M), \ \forall \alpha \in \Omega ^1(M) \ .
 \0(3.14)
$$

We recall that a symmetric linear connection $\nabla $ on a
Poisson manifold $(M,w)$ is called a {\it Poisson connection   }
if $\nabla w =0.$ Such connections exists iff $w$ is regular, i.e.
$rank \ w \ = \ const$ (see \cite{Va1}).

 PROPOSITION 3.9. {\it Let $(M,w)$ be a regular Poisson manifold
with a Poisson connection $\nabla .$ Then, the bivector $w^H,$
defined with respect to $\nabla ,$ is a Poisson structure on
$T^{*}M$ compatible with the canonical Poisson structure $W_0$ iff
the $2-$form
$$
(X,Y)\longrightarrow R_{\nabla }(X,Y)({\sharp}_w\alpha) \ , \ \
X,Y \in \chi (M)
 \0(3.15)
$$
vanishes for every Pfaff form $\alpha $ on $M.$ }

{\it Proof.} With $(3.5),$ the first condition $(3.12)$ becomes
$\nabla w=0,$ which we took as a hypotheses.
 The second condition $(3.12)$ becomes
$$
w^{ih}R^l_{hjk}=0 \ ,
$$
and we get the required conditions. \Endproof

 REMARK 3.10. If  $w$ is defined by a symplectic structure of $M$
then  $(3.15)$  means  $R_{\nabla }=0.$

\section*
{\begin{center}{\bf 4. Poisson structures derived from
differential forms}\end{center}}
 If $\omega $ is a $2-$form on
 a Riemannian manifold $(M,g)$ we
associate with it a $2-$form $\Theta (\omega )$ on the cotangent
bundle $\pi :T^{*}M\longrightarrow M,$ and considering  (pseudo-)
Riemannian metrics on $T^{*}M$ related to $g,$ we study the
conditions for  $\Theta (\omega )$   to produce  a Poisson
structure on this bundle.

Let $(M,g)$ be a $n-$dimensional manifold and $\nabla $ its
Levi-Civita connection. If $\Gamma ^k_{ij}$ are the local
coefficients of $\nabla ,$ a connection $\Gamma $ with the
coefficients $(3.5)$ is obtained on $T^{*}M.$

The system of local $1-$forms $(dx^i, \delta p_i),$ $(i=1,...,n),$
where
$$
\delta p_i:=dp_i+N_{ij}dx^j
 \0(4.1)
$$
defines the dual bases of the bases  $\{ {\delta}/{\delta x^i},
{\partial}/{\partial p_i} \}.$

The components of the curvature form are given by $(3.2).$ Since
the  connection is symmetric, the Bianchi identity $(3.13)$ holds.
The elements $\Phi ^k_{ij}$ of $(3.3)$ are
$$
\Phi ^k_{ij}=\Gamma ^k_{ij} \ .
 \0(4.2)
$$
The Riemannian metric $g$ provides  the "musical" isomorphism
${\sharp}_g:T^{*}M\longrightarrow $ $TM$
 and the codifferential
 $$
{\delta }_g:\Omega ^k(M)\longrightarrow \Omega ^{k-1}(M) \ , \  \
{({\delta}_g\alpha)}_{i_1...i_{k-1}}= -g^{st}{\nabla}_t
{\alpha}_{si_1...i_{k-1}} \ ,
    \0(4.3)
 $$
where $k\geq 1,$
$$
\alpha =\frac {1}{k!}{\alpha}_{i_1...i_k}dx^{i_1}\wedge ... \wedge dx^{i_k} \in
\Omega ^k(M) \
$$
and $(g^{st})$  are the entries of the inverse of the matrix
$(g_{ij})$  \cite{Va1}.

Let
$$
\omega =\frac 12 {\omega }_{ij}(x)dx^i \wedge dx^j \ , \ \ {\omega
}_{ji}=-{\omega }_{ij} \ ,
$$
be a $2-$form on $M.$

DEFINITION 4.1. The $2-$form $\Theta  (\omega )$ on $T^{*}M$ given
by
$$
\Theta  (\omega )={\pi}^{*}\omega -d\lambda \ ,
  \0(4.4)
$$
where $\lambda $ is the Liouville form, is said to be the {\it
associated $2-$form } of $\omega .$

With respect to the cobases $(dx^i, \delta p_i)$ we get
$$
\Theta (\omega ) = \frac 12 {\omega }_{ij}(x)dx^i \wedge dx^j +
dx^i \wedge \delta p_i \ .
   \0(4.5)
$$

Now, we   consider two (pseodo-) Riemannian metrics $G_1$ and
$G_2$ on $T^{*}M$ and study the conditions for the bivectors
$W_i={\sharp}_{G_i}\Theta (\omega ) ,$ $(i=1,2)$ to define
Poisson structures on $T^{*}M.$
 The Poisson condition $[W_i,W_i]=0,$ $i=1,2$ is equivalent to
 \cite{Va1}
 $$
{\delta}_{G_i}(\Theta (\omega ) \wedge \Theta (\omega ) )= 2\Theta
(\omega ) \wedge {\delta }_{G_i}\Theta (\omega ) \ ,
 \ \ i=1,2.
  \0(4.6)
 $$

 First,  consider \cite{{Vo1},{Vo2}}
  the pseudo-Riemannian metric $G_1$
 of signature $(n,n)$
$$
G_1=2\delta p_i\odot dx^i \ .
 \0(4.7)
$$
To find  the condition  which ensure that  $(4.6)$ holds, we need
the local expression of the codifferential ${\delta}_{G_1}$ of
$G_1.$ Denote by $\tilde{\nabla}$ the Levi-Civita connection  of
$G_1,$ and for simplicity we put
$$
{\tilde{\nabla}}_i={\tilde{\nabla}}_{\frac {\delta}{\delta x^i}} \
, \ \ {\tilde{\nabla}}^i={\tilde{\nabla}}_{\frac
{\partial}{\partial p_i}} \ .
  \0(4.8)
$$
$\tilde{\nabla}$ is defined by \cite{Vo1}
 $$
{\tilde{\nabla}}^i\frac {\partial}{\partial p_j}=0 \ , \ \
{\tilde{\nabla}}_i\frac {\partial}{\partial p_j}= -\Gamma
^j_{ik}\frac {\partial}{\partial p_k} \ ,
    \0(4.9)
 $$
$$
{\tilde{\nabla}}^i\frac {\delta}{\delta q^j}=0 \ , \ \
{\tilde{\nabla}}_i\frac {\delta}{\delta q^j}= \Gamma ^k_{ij}\frac
{\delta}{\delta q^k}- p_hR^h_{ijk}\frac {\partial}{\partial p_k} \
.
$$

 PROPOSITION 4.2. {\it The bivector ${\sharp}_{G_1}\Theta (\omega ) $
 defines a
Poisson structure on the cotangent bundle  $T^{*}M$ iff $\omega $
is a closed $2-$form on $M$ and $\Gamma ^a_{ai}=0,$ $\forall
i=1,...,n.$ In this case $\Theta (\omega ) $ is a symplectic form.
}

{\it Proof.} The proof is by a long computation in local
coordinates. After  computing the exterior product
 $\Theta (\omega ) \wedge \Theta (\omega )$ we get
 $$
{\delta}_{G_1}(\Theta (\omega ) \wedge \Theta (\omega ))= \frac
{2}{3!} \sum \limits _{(i,j,k)}{\nabla}_i{\omega}_{jk}dx^i\wedge
dx^j\wedge dx^k \ .
  \0(4.10)
 $$
Then, we  compute  ${\delta}_{G_1}\Theta (\omega ) $ and obtain
$$
\Theta (\omega ) \wedge {\delta}_{G_1}\Theta (\omega ) =\frac
{2}{3!} \sum \limits _{(i,j,k)} {\omega}_{ij}\Gamma
^a_{ak}dx^i\wedge dx^j\wedge dx^k+
    \0(4.11)
$$
$$
+({\delta}^k_j\Gamma ^a_{ai}-{\delta}^k_i\Gamma ^a_{aj})
dx^i\wedge dx^j \wedge \delta p_k \ .
$$
 $(4.6)$ implies
$$
{\delta}^k_j\Gamma ^a_{ai}-{\delta}^k_i\Gamma ^a_{aj}=0 \ , \ \
\forall i,j,k=1,...,n \ .
  \0(4.12)
$$
Making the contraction $k=j$ it follows that $\Gamma ^a_{ai}=0.$
Conversely, if $\Gamma ^a_{ai}=0$ then $(4.12)$ holds.
 Also, since $\nabla $ is symmetric, we get
$$
\sum \limits _{(i,j,k)}\frac {\partial {\omega }_{jk}}{\partial
x^i}= \sum \limits _{(i,j,k)} {\nabla}_i{\omega }_{jk} \ .
$$
Therefore, the condition $\sum \limits _{(i,j,k)}
 {\nabla}_i{\omega }_{jk} =0$ is equivalent to
 $d\omega =0.$
 \Endproof

Let us consider now the Riemannian metric of Sasaki type
$$
G_2=g_{ij}dx^i\odot dx^j+g^{ij}\delta p_i  \odot \delta p_j
  \0(4.13)
$$
(see \cite{Dom} for the Sasaki metric).

LEMMA 4.3. {\it The local coordinate expression of the Levi-Civita
connection $\bar{\nabla}$ of $G_2$ is
 $$
{\bar{\nabla}}^i\frac {\partial}{\partial p_j}=0 \ , \ \
{\bar{\nabla}}_i\frac {\partial}{\partial p_j}=- \frac
{1}{2}R^{jk}_{ \ \ i}\frac {\delta}{\delta q^k} -\Gamma
^j_{ik}\frac {\partial}{\partial p_k} \ ,
  \0(4.14)
 $$
$$
{\bar{\nabla}}^i\frac {\delta}{\delta q^j}= \frac 12R^{i \ \ k}_{
\ j}\frac {\delta}{\delta q^k}
 \ , \ \
{\bar{\nabla}}_i\frac {\delta}{\delta q^j}= \Gamma ^k_{ij}\frac
{\delta}{\delta q^k}-\frac 12 R_{kij}\frac {\partial}{\partial
p_k} \ ,
$$
where we used again the notations of $(4.8)$ and $R^{jk}_{ \ \ i}$
(also $R^{i \ \ k}_{ \ j}$) are obtained from $R_{kij}$ by the
operation of lifting the indices, i.e.}
$$
R^{jk}_{ \ \ i}=g^{ja}g^{kb}R_{abi} \ , \ \ R^{i \ \ k}_{ \ j}=
g^{ia}g^{kb}R_{ajb} \ .
$$

{\it Proof.} The result is proved by a straightforward
computation. \Endproof

 PROPOSITION 4.4. {\it The bivector ${\delta}_{G_2}\Theta (\omega ) $
defines a Poisson structure on the cotangent bundle $T^{*}M$ iff
$$
\nabla \omega =0 \ , \ \ g^{ab}R^k_{abi}=0 \ , \ \
{\omega }^{ab}R^k_{iab}=0 \ ,
   \0(4.15)
$$
where ${\omega}^{ab}=g^{ai}g^{bj}{\omega }_{ij}$ are the
components of the bivector
 $w={\sharp}_g\omega $ on $M.$}

 {\it Proof:} By a new long computation again, we get
$$
\frac 12{\delta}_{G_2}(\Theta (\omega ) \wedge \Theta (\omega ))=
\frac {1}{3!} g^{ab}{\nabla }_a(\sum \limits _{(i,j,k)}{\omega
}_{ij} {\omega }_{kb})dx^i\wedge dx^j\wedge dx^k-
$$
$$
-(g^{ab}\sum \limits _{(i,j,k)}({\nabla }_a{\omega }_{ij}{\delta
}^k_b) dx^i\wedge dx^j \wedge \delta p_k+ \frac 12{\omega
}_{ab}(R^{kab}{\delta }^j_i-R^{jab}\delta ^k_i) dx^i \wedge \delta
p_j \wedge \delta p_k
$$
and
$$
\Theta (\omega ) \wedge {\delta}_{G_2}\Theta (\omega )= \frac
{1}{3!}\sum \limits _{(i,j,k)}{({\delta}_{G_2}\Theta (\omega ))}_k
dx^i\wedge dx^j\wedge dx^k+
$$
$$
+\frac {1}{2!}[\delta ^k_i{({\delta}_{G_2}\Theta (\omega ))}_j-
\delta ^k_j{({\delta}_{G_2}\Theta (\omega ))}_i] dx^i\wedge
dx^j\wedge \delta p_k \ ,
$$
where
$$
{{\delta}_{G_2}\Theta (\omega )}={({\delta}_{G_2}\Theta (\omega
))}_kdx^k= g^{ab}({\nabla }_a{\omega }_{kb}-\frac 12R_{abk})dx^k \
.
$$
Identifying the coefficients, the Poisson conditions $(4.6)$ for
$W_2$ becomes:
$$
 g^{ab}\sum \limits _{(i,j,k)}{\omega }_{ij}R^h_{abk}=0  \ , \ \
g^{ab}\sum \limits _{(i,j,k)}({\nabla }_a{\omega }_{ij})\omega
_{kb}=0 \ ,
  \0(4.16)
$$
$$
 \nabla \omega =0 \ , \ \  \ g^{ab}R^k_{abi}=0 \ ,
  \0(4.17)
$$
and
$$
\omega ^{ab}R^k_{iab}=0 \ .
  \0(4.18)
$$
Let us remark that the conditions $(4.17)$ implies $(4.16),$
 because, if
 $\nabla \omega =0$ then
 ${\nabla }_a\omega _{ij}=0,$ and  $g^{ab}R^k_{abi}=0$ implies
 $g^{ab}{\omega }_{ij}R^h_{abk}=0.$
 \Endproof

 REMARK 4.5. If the bivector ${\sharp }_{G_2}\Theta (\omega )$ defines a
Poisson structure on $T^{*}M$ then $w={\sharp }_g\omega $ defines
a Poisson structure on $M,$ as the second  condition $(4.16)$ is
equivalent to the Poisson condition \cite{Va1}
$$
\sum \limits _{(i,j,k)}w^{ia}{\nabla}_aw^{jk}=0 \ .
$$
(The local coordinate expression of $w$ is $(1.5).$)

 COROLLARY 4.6. {\it
If ${\sharp }_{G_2}\Theta (\omega )$ is a Poisson bivector
 on $T^{*}M,$ then the scalar curvature $r$ of $(M,g)$ vanishes.}

 {\it Proof: }
The expression of $r$ is $ r=g^{ab}R_{ab} \ , $ where
$R_{ba}=R^k_{akb}=R_{ab}$ are the components of the Ricci tensor,
and
 if we make the contraction $k=i$ in the second relation $(4.15)$ we
 get
 $g^{ab}R^k_{akb}=0,$ and whence $r=0.$
 \Endproof

\vspace{1mm} {\bf Acknowledgement} During the work on this paper,
the author  held a postdoctoral grant at the University of Haifa,
Israel. He would like to thank the University of Haifa, its
Department of Mathematics, and, personally, Prof. Izu Vaisman for
suggestions  and hospitality during the period when the
postdoctoral program was completed.

{\small \begin{tabular}{ll}
Gabriel Mitric   \\
Catedra de Geometrie\\
Universitatea ``Al. I. Cuza" \\
Ia\c{s}i 6600
\\Rom\^{a}nia\\
gmitric@uaic.ro
\end{tabular}}

\end{document}